\documentclass[11pt]{article}

\usepackage{amsmath,amssymb,amsthm,amsfonts,verbatim}
\usepackage{graphicx}
\usepackage{enumerate}
\usepackage[all,2cell]{xy}

\newtheorem{theorem}{Theorem}[section]

\newtheorem{problem}[theorem]{Problem}
\theoremstyle{definition}
\newtheorem{example}[theorem]{Example}

\theoremstyle{definition}
\newtheorem{definition}[theorem]{Definition}
\newtheorem{remark}[theorem]{Remark}

\newcommand{\nc}{\newcommand}
\nc{\dmo}{\DeclareMathOperator}

\nc{\FF}{\mathbf{F}}
\nc{\Q}{\mathbb{Q}}
\nc{\R}{\mathbb{R}}
\nc{\Z}{\mathbb{Z}}
\nc{\C}{\mathbb{C}}
\nc{\N}{\mathbb{N}}
\nc{\I}{\mathcal{I}}
\nc{\bwedge}{\textstyle{\bigwedge}}
\nc{\p}{{\frak p}}

\nc{\astT}[1]{\ensuremath{\ast_{#1}}}
\def\Fqbar{\overline{\F}_q}
\def\FF{{\cal F}}
\dmo{\tr}{tr}
\dmo\roots{Roots}
\dmo{\GL}{GL}
\dmo{\Sp}{Sp}
\dmo{\Aut}{Aut}
\dmo{\Stab}{Stab}
\dmo\SL{SL}

\dmo\im{im}
\dmo\id{id}
\dmo\Sym{Sym}
\dmo\End{End}
\dmo\Conf{Conf}
\dmo\UConf{UConf}
\dmo\op{op}
\dmo\coker{coker}
\dmo\Inj{Inj}
\dmo\Map{Map}
\dmo\chr{char}
\dmo\Ext{Ext}
\dmo\spn{span}
\dmo\Ind{Ind}
\dmo\Res{Res}
\dmo\Out{Out}
\dmo\Mod{Mod}
\dmo\PMod{PMod}
\dmo\Fix{Fix}
\dmo\Frob{Frob}

\dmo\F{\mathbb{F}}
\dmo\IA{IA}

\def\fsl{\mathfrak{sl}}

\def\fsp{\mathfrak{sp}}

\dmo\gr{gr}

\dmo{\Rep}{Rep}

\nc{\cR}{\mathcal{R}}
\nc{\cM}{\mathcal{M}}
\nc{\cC}{\mathcal{C}}

\nc{\ra}{\rightarrow}
\nc{\inj}{\hookrightarrow}
\nc{\surj}{\twoheadrightarrow}

\nc{\margin}[1]{\marginpar{\scriptsize #1}}
\nc{\para}[1]{\medskip\noindent\textbf{#1.}}

\title{Representation stability}

\author{Benson Farb\thanks{The author gratefully acknowledges support from the National Science
    Foundation. He would also like to thank his collaborators Thomas Church and Jordan Ellenberg, 
    without whom the work discussed herein would not exist.  Thanks also to Thomas Church and Dan Margalit for their extensive comments on an earlier version of this paper.}}

\begin{document}

\maketitle

\begin{abstract}
Representation stability is a phenomenon whereby the structure of certain sequences $X_n$ of spaces can be seen to stabilize when viewed through the lens of representation theory.  In this paper I  describe this phenomenon and sketch a framework, the theory of FI-modules, that explains the mechanism behind it.   
\end{abstract}

\section{Introduction}

Sequences $V_n$ of representations of the symmetric group $S_n$ occur naturally in topology, combinatorics, algebraic geometry and elsewhere.  Examples include the cohomology of configuration spaces $\Conf_n(M)$, moduli spaces of $n$-pointed Riemann surfaces and 
congruence subgroups $\Gamma_n(p)$; spaces of polynomials on rank varieties of $n\times n$ matrices; and $n$-variable diagonal co-invariant algebras.   

Any $S_n$-representation is a direct sum of irreducible representations.  These are parameterized by partitions of $n$.   Following a 1938 paper of Murnaghan, one can pad a partition 
$\lambda=\sum_{i=1}^r\lambda_i$ of any number $d$ to produce a partition $(n-|\lambda|)+\lambda$ of $n$ for all $n\geq |\lambda|+\lambda_1$.  The decomposition of $V_n$ into irreducibles thus produces a sequence of multiplicities 
of partitions $\lambda$, recording how often $\lambda$ appears in $V_n$.

A few years ago Thomas Church and I discovered that for many important sequences $V_n$ arising 
in topology, these multiplicities become constant once $n$ is large enough.   With Jordan Ellenberg and Rohit Nagpal, we built a theory to explain this stability, converting  it to a finite generation property for a single object.    We applied this to prove stability in these and many other examples.   As a consequence, the character of $V_n$ is given (for all $n\gg 1$) by a single polynomial, called a {\em character polynomial}, studied  by Frobenius but not so widely known 
today.   One of the main points of our work is that the mechanism for this stability comes from a common structure underlying all of these examples.  

After giving an overview of this theory, we explain how it applies and 
connects to an array of counting problems for polynomials over finite fields, and for maximal tori in the finite groups $\GL_n\F_q$.  In particular, the stability of such counts reflects, and is reflected in, the representation stability of the cohomology of an associated algebraic variety.  We begin with a motivating example.

\section{Configuration spaces and representation theory}
\label{section:configs}

Let $M$ be any connected, oriented manifold.  For any $n\geq 1$ 
let $\Conf_n(M)$ be the space of configurations of ordered $n$-tuples of distinct points in $M$:
\[\Conf_n(M):=\{(z_1,\ldots ,z_n)\in M^n :  z_i\neq z_j \ \text{if $i\neq j$}\}.\]

The symmetric group $S_n$ acts freely on $\Conf_n(M)$ by permuting the coordinates:
\[\sigma\cdot (z_1,\ldots ,z_n):=(z_{\sigma(1)},\ldots ,z_{\sigma(n)}).\]
This action induces for each $i\geq 0$ an action of $S_n$ on the complex vector space $H^i(\Conf_n(M);\C)$, making   $H^i(\Conf_n(M);\C)$ into an $S_n$-representation.
Here we have chosen $\C$ coefficients for simplicity of exposition.

\subsection{Cohomology of configuration spaces} The study of  configuration spaces and their cohomology is a classical topic. We concentrate on the following fundamental problem. 

\begin{problem}[{\bf Cohomology of configuration spaces}]
\label{problem:main1}
Let $M$ be a connected, oriented manifold. Fix a ring $R$. Compute $H^*(\Conf_n(M);R)$ as an $S_n$-representation.
\end{problem}

Problem~\ref{problem:main1} was considered in special cases by Brieskorn, F. Cohen, Stanley, Orlik, Lehrer-Solomon and many others; see, e.g.\ \cite{LS} and the references contained therein.  

What exactly does ``compute as an $S_n$-representation'' mean?   Well, by Maschke's Theorem,  any $S_n$-representation over $\C$ is a direct sum of irreducible $S_n$-representations.   In 1900, Alfred Young gave an explicit bijection between the set of (isomorphism classes of) irreducible $S_n$-representations  and the set of partitions $n=n_1+\cdots +n_r$ of $n$ with $n_1\geq \cdots \geq n_r>0$.  

Let $\lambda=(a_1,\ldots ,a_r)$ be an $r$-tuple of integers with $a_1\geq \cdots \geq a_r>0$ and such that $n-\sum a_i\geq a_1$.  We denote by $V(\lambda)$, or $V(\lambda)_n$ when we want to emphasize $n$, the representation in Young's classification corresponding to the partition 
$n=(n-\sum_{i=1}^ra_i) +a_1+\cdots +a_r$.  With this terminology we have, for 
example, that $V(0)$ is the trivial representation, $V(1)$ is the $(n-1)$-dimensional irreducible representation 
$\{(z_1,\ldots ,z_n)\in\C^n:\sum z_i=0\}$, and $V(1,1,1)$ is the irreducible representation $\bigwedge^3V(1)$.   For each $n\geq 1$ and each $i\geq 0$ we can write
\begin{equation}
\label{eq:generalsum1}
H^i(\Conf_n(M);\C)=\bigoplus_{\lambda}d_{i,n}(\lambda) V(\lambda)_n
\end{equation}
for some integers $d_{i,n}(\lambda)\geq 0$.  The sum on the right-hand side of \eqref{eq:generalsum1} is taken over all partitions $\lambda$ of numbers $\leq n$ for which $V(\lambda)_n$ is defined.  The coefficient $d_{i,n}(\lambda)$ is called the {\em multiplicity} of $V(\lambda)$ in $H^i(\Conf_n(M);\C)$.   

\bigskip
\noindent
{\bf Problem~\ref{problem:main1} over $\boldmath\C$, restated: } Compute the multiplicities $d_{i,n}(\lambda)$.

\bigskip
Why should we care about solving this problem?  Here are a few reasons:

\begin{enumerate}
\item Even the multiplicity $d_{i,n}(0)$ of the trivial representation $V(0)$ is interesting: it computes the 
$i^{\rm th}$ Betti number of the space $\UConf_n(M):=\Conf_n(M)/S_n$ of unordered $n$-tuples of distinct points in $M$.  In other words, 
\begin{equation}
\label{eq:transfer1}
\dim_\C H^i(\UConf_n(M);\C)=d_{i,n}(0)
\end{equation}
by transfer applied to the finite cover $\Conf_n(M)\to \UConf_n(M)$.

\item More generally, the $d_{i,n}(\lambda)$ for other partitions $\lambda$ of $n$ compute the Betti numbers of other (un)labelled configuration spaces.  For example, for fixed $a,b,c\geq 0$, consider the space $\Conf_n(M)[a,b,c]$ of configurations of $n$ distinct labelled points on $M$ where one colors $a$ of the points blue, $b$ red, and $c$ yellow, and where points of the same color are indistinguishable from each other.  Then $H_i(\Conf_n(M)[a,b,c];\C)$ can be determined from $d_{i,n}(\mu)$ for certain $\mu=\mu(a,b,c)$.   See \cite{Ch} for a discussion, and see \cite{VW} for an explanation of how these 
spaces arise naturally in algebraic geometry. 

\item The representation theory of $S_n$ provides strong constraints on the possible values of 
$\dim_\C H^i(\Conf_n(M);\C)$.  As a simple example, if the action of $S_n$ on $H^i(\Conf_n(M);\C)$ is {\em essential} in a specific sense (cf.\ \S\ref{section:intro:instability}) for $n\gg 1$, then one can conclude for purely representation-theoretic reasons that $\lim_{n\to\infty}\dim_\C H^i(\Conf_n(M);\C)=\infty$.   
This happens for example for every $i\geq 1$ when $M=\C$.  See \S\ref{section:intro:instability} below.

\item For certain special smooth projective varieties $M$, the multiplicities $d_{i,n}(\lambda)$ encode and are encoded by  delicate  information about the combinatorial statistics of the 
$\F_q$-points of $M$ or related varieties; see \S\ref{section:varieties} below for two specific applications.

\item The decomposition \eqref{eq:generalsum1} can have geometric meaning, and can point the way for us to guess at meaningful topological invariants.  We now discuss this in a specific example.

\end{enumerate}

\subsection{A case study: the invariants of loops of configurations}  Consider the special case where $M$ is the complex plane $\C$.  Elements of $H^1(\Conf_n(\C);\C)$ are homomorphisms $\pi_1(\Conf_n(\C))\to\C$.  Computing $H^1(\Conf_n(\C);\C)$ is thus answering the basic question: 

\medskip
{\it What are the ways of attaching a complex number to each loop of configurations of $n$ points in the plane, in a way that is natural (= additive)?}

\begin{figure}[h]
    \centering
      \includegraphics[width=80mm]{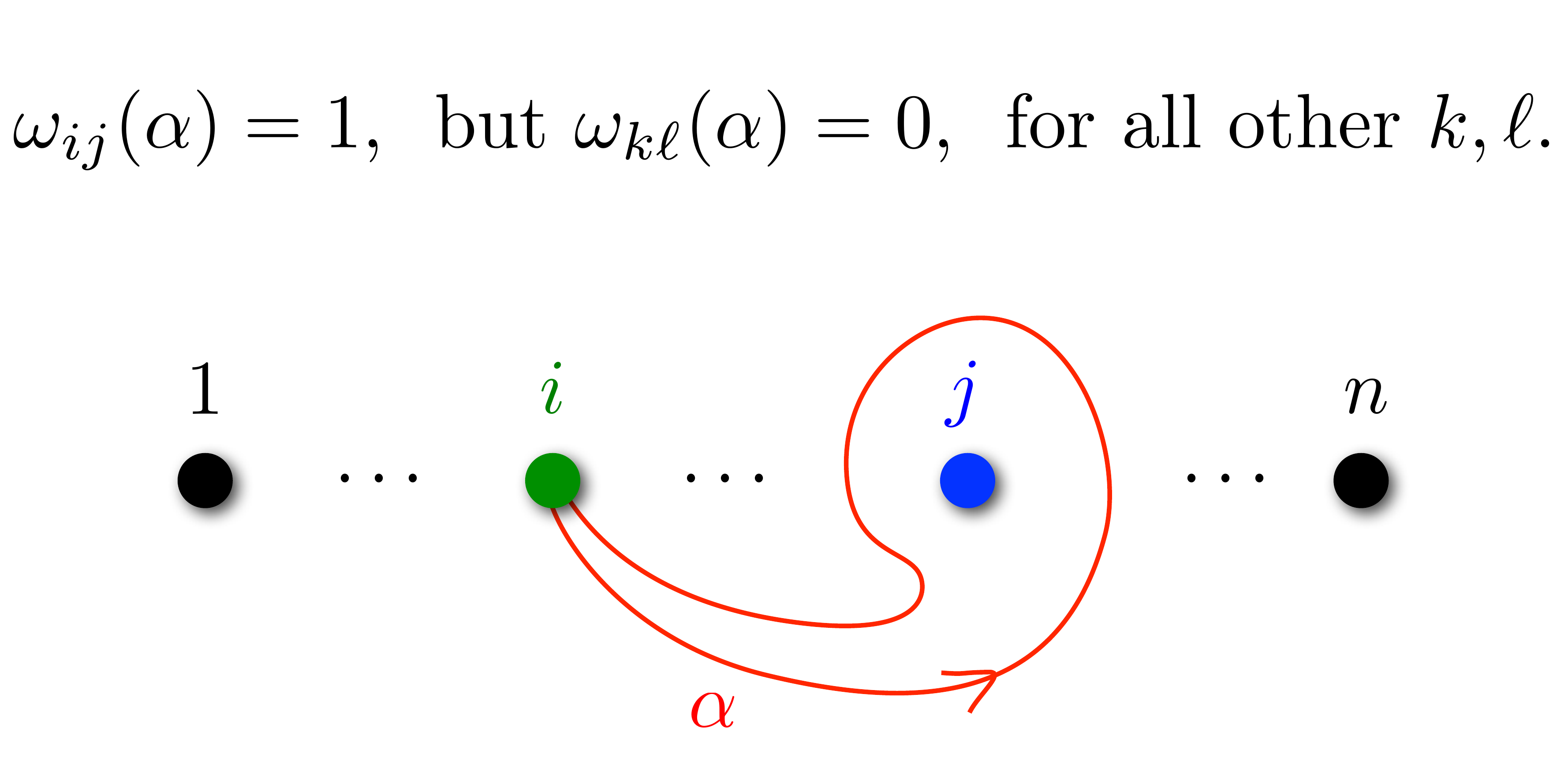}
    \caption{The proof that $\{\omega_{ij}: i<j\}$ is linearly independent in $H^1(\Conf_n(\C);\C)$.}
     \label{fig:braid1}
  \end{figure}

\medskip
To construct examples, pick $1\leq i,j \leq n$ with $i\neq j$.  Given any loop $\gamma(t)=(z_1(t),\ldots ,z_n(t))$ in $\Conf_n(\C)$, we can ignore all points except for $z_i(t)$ and $z_j(t)$ and measure 
how much $z_j(t)$ winds around $z_i(t)$; namely we let $\alpha_{ij}:[0,1]\to \C$ be the loop $\alpha_{ij}(t):=z_j(t)-z_i(t)$ and set
\[\omega_{ij}(\gamma):=\frac{1}{2\pi i}\int_{\displaystyle \alpha_{ij}}\frac{dz}{z}\]
It is easy to verify that $\omega_{ij}:\pi_1(\Conf_n(\C))\to\C$ is indeed a homomorphism, that $\omega_{ij}=\omega_{ji}$, and that the set $\{\omega_{ij}: i<j\}$ is linearly independent in $H^1(\Conf_n(\C);\C)$; see Figure~\ref{fig:braid1}.  
Linear combinations of the $\omega_{ij}$ are in fact the only natural invariants of loops of configurations in $\C$.

\begin{theorem}[{\bf Artin(1925), Arnol'd(1968) \cite{Ar}}]
The set $\{\omega_{ij}: 1\leq i<j\leq n\}$ is a basis for $H^1(\Conf_n(\C);\C)$ for any $n\geq 2$.  Thus  $H^1(\Conf_n(\C);\C)\approx \C^{\binom{n}{2}}$.
\end{theorem}

There is more to say.   The $S_n$ action on $H^1(\Conf_n(\C);\C)$ is determined by its action on the basis via $\sigma\cdot \omega_{ij}=\omega_{\sigma(i)\sigma(j)}$, from which we can deduce that

\begin{equation}
\label{eq:Artin}
H^1(\Conf_n(\C);\C)=V(0)\oplus V(1)\oplus V(2) \ \ \ \ \text{for $n\geq 4$}
\end{equation}
 using only elementary representation theory.  We can see from this algebraic picture that the subspace  of vectors fixed by all of $S_n$ is $1$-dimensional, spanned by the 
 vector \[\Omega:=\sum_{1\leq i<j\leq n}\omega_{ij}\in H^1(\Conf_n(\C);\C).\]   
 This implies the following geometric statement:  the only 
 natural invariant of loops of configurations of $n$ distinct {\em unordered} points in $\C$ is total winding number $\Omega$; in particular, $H^1(\Conf_n(\C)/S_n;\C)\approx\C$.   
 
Looking again at \eqref{eq:Artin} we see a copy of the standard permutation representation $\C^n=V(0)\oplus V(1)$ given by $\sigma\cdot u_i=u_{\sigma(i)}$, with $u_i=\sum_{j\neq i}\omega_{ij}$.   This indicates that the $u_i$ should be geometrically meaningful, which indeed they are: $u_i$ gives 
the total winding number of all points $z_j$ around the point $z_i$.   

I hope that even the simple example of $\Conf_n(\C)$ convinces the reader that understanding $H^i(\Conf_n(M);\C)$ as an $S_n$-representation and not just as a naked vector space gives us a much richer geometric picture. 

\subsection{Homological (in)stability}  
\label{section:intro:instability}
The above discussion fits in to a broader context.  Let $X_n$ be a sequence of spaces or groups.  The classical theory of (co)homological stability (over a fixed  ring $R$) in topology produces results of the form: the homology $H_i(X_n;R)$ (resp.\ $H^i(X_n;R)$) does not depend on $n$ for $n\gg i$.   This converts an {\it a priori} infinite computation to a finite one.  Examples of such sequences $X_n$ include symmetric groups $S_n$ (Nakaoka), braid groups $B_n$ (Arnol'd and F. Cohen), the space $\UConf_n(M)$ of $n$-point subsets of the interior of a compact, connected manifold $M$ with nonempty boundary (McDuff, Segal), special linear groups $\SL_n\Z$ (Borel and Charney), the moduli space ${\cal M}_n$ of genus $n\geq 2$ Riemann surfaces (Harer), and automorphism groups of free groups (Hatcher-Vogtmann-Wahl) ; see \cite{Co} for a survey. 

For many natural sequences $X_n$ homological stability fails in a strong way. We saw above that $H^1(\Conf_n(\C);\C)\approx \C^{\binom{n}{2}}$ for all $n\geq 4$.   In fact, one can prove for each $i\geq 1$  that 
\begin{equation}
\label{eq:instability1}
\lim_{n\to\infty}\dim_\C(H^i(\Conf_n(\C);\C))=\infty .
\end{equation}
The underlying mechanism behind this instability is symmetry.  Call a representation $V$ of $S_n$ {\em not essential}  if  $\sigma\cdot v=\pm v$ for each $\sigma\in S_n$ and each $v\in V$.   Basic representation theory of $S_n$ implies that any essential representation $V$ of $S_n,n\geq 5$ satisfies $\dim(V)\geq n-1$.   It is not hard to check that the $S_n$-representation  $H^i(\Conf_n(\C);\C))$ is essential, implying the blowup \eqref{eq:instability1}.    One hardly needs to know topology to prove \eqref{eq:instability1}! The driving force behind this is the representation theory of $S_n$.

More generally, whenever we have a sequence of larger and larger groups  $G_n$, and a sequence $V_n$ of ``essential'' $G_n$-representations, one expects that $\dim(V_n)\to\infty$.   The general slogan of representation stability is:  in many situations the {\em names} of the representations $V_n$ should stabilize as $n\to\infty$.  The question is, how can we formalize this slogan, and how can we use such information?  We focus on the case $G_n=S_n$; see \S\ref{section:other:examples} for a discussion of other examples, such 
as  $G_n=\SL_n\F_p, \GL_n\Z$ and $\Sp_{2n}\Z$.  


\section{Representation stability (the $S_n$ case)} 
\label{section:repstab}
With our notation, the description of $H^1(\Conf_n(\C);\C)$ given in \eqref{eq:Artin} does not depend on $n$ once $n\geq 4$.    In 2010 Thomas Church and I guessed that such a phenomenon might be true for cohomology in all degrees.  Using an inductive description of $H^\ast(\Conf_n(\C);\C)$ as a sum of induced representations (a weak form of a theorem of Lehrer-Solomon \cite{LS}), one can convert this question into a purely representation-theoretic one.   After doing this, we asked David Hemmer about the case $i=2$.  He wrote a computer program that produced the following output; we use the notation $C_n$ for $\Conf_n(\C)$ to save space.
\begin{equation}
\label{eq:stabilizationinaction}
\begin{array}{l}
  H^2(C_4;\C) = V(1)^{\oplus 2}\oplus V(1,1) \oplus V(2)\\  
  \\
  H^2(C_5;\C) = V(1)^{\oplus 2} \oplus V(1,1)^{\oplus 2}
  \oplus V(2)^{\oplus 2} \oplus 
  V(2,1)\\  
  \\
  H^2(C_6;\C) = V(1)^{\oplus 2} \oplus  V(1,1)^{\oplus 2}
  \oplus V(2)^{\oplus 2} \oplus 
  V(2,1)^{\oplus 2} \oplus V(3)\\ 
  \\
  H^2(C_7;\C) = V(1)^{\oplus 2} \oplus V(1,1)^{\oplus 2}
  \oplus V(2)^{\oplus 2} \oplus 
  V(2,1)^{\oplus 2} \oplus  V(3) \oplus V(3,1)\\
  \\
  H^2(C_8;\C) = V(1)^{\oplus 2} \oplus V(1,1)^{\oplus 2}
  \oplus V(2)^{\oplus 2} \oplus 
  V(2,1)^{\oplus 2} \oplus  V(3) \oplus V(3,1)\\
  \\
  H^2(C_9;\C) = V(1)^{\oplus 2} \oplus V(1,1)^{\oplus 2}
  \oplus V(2)^{\oplus 2} \oplus 
  V(2,1)^{\oplus 2} \oplus  V(3) \oplus V(3,1)\\
  \\
H^2(C_{10};\C) = V(1)^{\oplus 2} \oplus V(1,1)^{\oplus 2}
  \oplus V(2)^{\oplus 2} \oplus 
  V(2,1)^{\oplus 2} \oplus  V(3) \oplus V(3,1)
\end{array}
\end{equation}

This was compelling.   Indeed, it turns out that this decomposition holds for $H^2(C_n;\C)$ for all $n\geq 7$, so  the decomposition of $H^2(\Conf_n(\C);\C)$ into irreducible $S_n$-representations {\em stabilizes}.  The most useful way we found to encode this type of behavior was via the notion of a representation stable sequence, which we now explain.

\subsection{Representation stability and $H^i(\Conf_n(M);\C)$}
Let $V_n$ be a sequence of $S_n$--representations equipped with
linear maps $\phi_n\colon V_n\to V_{n+1}$ so that following diagram commutes for each $g\in S_n$:
\[\xymatrix{
  V_n\ar^{\phi_n}[r]\ar_{g}[d]&V_{n+1}\ar^{g}[d]\\
  V_n\ar_{\phi_n}[r]&V_{n+1} }\]
Here $g$ acts on $V_{n+1}$ by its image under the standard inclusion
$S_n\hookrightarrow S_{n+1}$.  We call such a sequence of
representations \emph{consistent}. We made the following definition in \cite{CF1}.

\begin{definition}[{\bf Representation stability for $S_n$-representations}] 
\label{definition:repstab}
  A consistent sequence $\{V_n\}$ of $S_n$--representations is \emph{representation stable} if there exists $N>0$ so that for all $n\geq N$, each of the following conditions holds:
\begin{enumerate}
\item \textbf{Injectivity:} The maps $\phi_n\colon V_n\to V_{n+1}$ are 
  injective.
\item \textbf{Surjectivity:} The span of the $S_{n+1}$--orbit of
  $\phi_n(V_n)$ is all of $V_{n+1}$.
\item \textbf{Multiplicities:} Decompose $V_n$ into irreducible
  $S_n$--representations as
  \[V_n=\bigoplus_\lambda c_n(\lambda)V(\lambda)_n\] with
  multiplicities $0\leq c_n(\lambda)\leq \infty$. Then $c_n(\lambda)$ does not depend on $n$.  
  \end{enumerate}
\end{definition}

The number $N$ is called the {\em stable range}.  The sequence $V_n:=\wedge^*\C^n$ of exterior algebras is an example of a consistent sequence of $S^n$-representations that is not representation stable.

It is not hard to check that, given Condition 1 for $\phi_n$, Condition 2 for $\phi_n$ is equivalent to the following : $\phi_n$ is a composition of the inclusion $V_n\hookrightarrow \Ind_{S_n}^{S_{n+1}}V_n$ with a surjective $S_{n+1}$--module homomorphism $ \Ind_{S_n}^{S_{n+1}}V_n\to V_{n+1}$.   This point of view leads to the stronger condition {\em central stability}, a very useful concept invented by Putman \cite{Pu} at around the same time, which he applied in his study of the cohomology of congruence subgroups. 

There are variations on Definition~\ref{definition:repstab}.  For example one can allow the 
stable range $N$ to depend on the partition $\lambda$.  In \cite{CF1} we define representation stability for other sequences $G_n$ of groups, with a definition analogous to Definition~\ref{definition:repstab} with $G_n=S_n$ replaced by $G_n=\GL_n\Z$, $\Sp_{2g}\Z, \GL_n\F_q, \Sp_{2g}\F_q$, and 
hyperoctahedral groups $W_n$; see \S\ref{section:other:examples} below.  In each case one needs a coherent naming system for representations of $G_n$ as $n$ varies.  

\begin{remark}
We originally stumbled onto representation stability in \cite{CF2} while making some computations in the homology of the Torelli group $\I_g$.  In this situation the homology $H_i(\I_g;\C)$ is a representation of the integral symplectic group $\Sp_{2g}\Z$.  We found some $\Sp_{2g}\Z$-submodules of $H_i(\I_g;\C)$ whose names did not depend on $g$ for $g\gg 1$.  Representation stability (for sequences of $\Sp_{2g}\Z$-representations) arose from our attempt to formalize this.  See \cite{CF2}.   After \cite{CF1,CF2} appeared, Richard Hain kindly shared with us some of his unpublished notes from the early 1990s, where he also developed a conjectural picture of the homology $H_i(\I_g;\C)$ as an $\Sp_{2g}\Z$-representation that is similar to the idea of representation stability for $\Sp_{2g}\Z$-representations presented in \cite{CF1}.
\end{remark}

Using in a crucial way a result of Hemmer \cite{He}, we proved in \cite{CF1}  the following.
\begin{theorem}[{\bf Representation stability for \boldmath$\Conf_n(\C)$}]
\label{theorem:CF1}
For any fixed $i\geq 0$, the sequence $\{H^i(\Conf_n(\C);\C)\}$ is representation stable with stable range $n\geq 4i$.   
\end{theorem}
The stable range $n\geq 4i$ given in Theorem~\ref{theorem:CF1} predicts that $H^2(\Conf_n(\C);\C)$ will stabilize once $n=8$; in truth it stabilizes starting at $n=7$.

The problem of computing all of the stable multiplicities $d_{i,n}(\lambda)$ in the decomposition of $H^i(\Conf_n(\C);\C)$ is thus converted to a problem which is finite and in principle solvable by a computer.  However, putting this into practice is a delicate matter, and the actual answers can be quite complicated.  For example, for $n\geq 16$: 
\begin{align*}
H&{}^4(\Conf_n(\C);\C)=\quad\\
\ &V(1)^{\oplus 2}
\oplus V(2)^{\oplus 6}
\oplus V(1,1)^{\oplus 6}
\oplus V(3)^{\oplus 8}
\oplus V(1,1,1)^{\oplus 9}
\oplus V(2,1)^{\oplus 16}\\
\oplus\,&V(4)^{\oplus 6}
\oplus V(1,1,1,1)^{\oplus 5}
\oplus V(5)^{\oplus 2}
\oplus V(2,2)^{\oplus 12}
\oplus V(3,1)^{\oplus 19}\\
\oplus \,&V(2,1,1)^{\oplus 17}
\oplus V(4,1)^{\oplus 12}
\oplus V(2,1,1,1)^{\oplus 7}
\oplus V(3,2)^{\oplus 14}
\oplus V(2,2,1)^{\oplus 10}\\
\oplus \,&V(5,1)^{\oplus 3}
\oplus V(3,3)^{\oplus 4}
\oplus V(3,1,1)^{\oplus 16}
\oplus V(2,2,2)^{\oplus 2}
\oplus V(4,2)^{\oplus 7}\\
\oplus \,&V(4,1,1)^{\oplus 8}
\oplus V(5,2)
\oplus V(2,2,1,1)^{\oplus 2}
\oplus V(3,1,1,1)^{\oplus 5}
\oplus V(5,1,1)^{\oplus 2}\\
\oplus \,&V(4,3)^{\oplus 2}
\oplus V(3,2,1)^{\oplus 9}
\oplus V(4,1,1,1)^{\oplus 2}
\oplus V(3,3,1)^{\oplus 2}
\oplus V(3,2,2)\\
\oplus \,&V(4,2,1)^{\oplus 3}
\oplus V(3,2,1,1)
\oplus V(5,1,1,1)
\oplus V(4,3,1)
\end{align*}

Theorem~\ref{theorem:CF1} was greatly extended by Church in \cite{Ch} from $M=\C$ to $M$ any connected, oriented manifold, as follows.  

\begin{theorem}[{\bf Representation stability for configuration spaces}]
\label{theorem:church}
Let $M$ be any connected, oriented manifold with $\dim(M)\geq 2$ and with $\dim_\C H^*(M;\C)<\infty$.  Fix $i\geq 0$. Then the sequence $\{H^i(\Conf_n(M);\C)\}$ is representation stable with stable range  $n\geq 2i$ if $\dim(M)\geq 3$ and $n\geq 4i$ if $\dim(M)=2$.
\end{theorem}

This of course leaves open the following.

\begin{problem}[{\bf Computing stable multiplicities}]
Given a connected, oriented manifold $M$, compute explicitly the stable multiplicities 
$d_{i,n}(\lambda)$ for the decomposition of $H^i(\Conf_n(M);\C)$ into irreducibles.    Give geometric interpretations of these numbers, as in the case of $H^1(\Conf_n(\C);\C)$ discussed in \S\ref{section:configs} above.
\end{problem}

The problem of computing the $d_{i,n}(\lambda)$ seems to have been solved in very few cases. For example, I do not know the answer even for $M$ a closed surface of genus $g\geq 1$.   

The paper \cite{CF1} gives many other examples of representation stable sequences $V_n$ that arise naturally in mathematics, from the cohomology of Schubert varieties to composition of Schur functors to many of the examples given in \S\ref{section:fiexamples} below.

\subsection{An application to classical homological stability} 
Consider the space $\UConf_n(M):=\Conf_n(M)/S_n$ of unordered 
$n$-tuples of distinct points in $M$. As mentioned above, when $M$ is the interior of a compact manifold with nonempty boundary, classical homological stability  for $H_i(\UConf_n(M);\Z)$ was proved by McDuff and Segal, generalizing earlier work of Arnol'd and F. Cohen.   The reason that 
the assumption $\partial M\neq \emptyset$ is needed is that in this case one has a map $\psi_n:\Conf_n(M)\to\Conf_{n+1}(M)$ for each $n\geq 0$ given by ``injecting a point at infinity'' (see Proposition 4.6 of \cite{CEF1} for details).  While $\psi_n$ is really just defined up to homotopy, it induces for each $i\geq 0$ a well-defined homomorphism 
\begin{equation}
\label{eq:mcduff}
(\psi_n)_*:H_i(\UConf_n(M);\Z)\to H_i(\UConf_{n+1}(M);\Z)
\end{equation}
which McDuff and Segal prove is an isomorphism for $n\geq 2i+2$.  This is the typical way one proves classical homological stability for a sequence of spaces $X_n$, namely one finds maps $X_n\to X_{n+1}$ and proves that they eventually induce isomorphisms on homology.  

What about the case when $M$ is closed? In this case there are no natural maps between $\UConf_n(M)$ and $\UConf_{n+1}(M)$.  A natural thing to do would be to consider the $S_n$-cover $\Conf_n(M)\to \UConf_n(M)$, where there are maps (in fact $n+1$ of them) $\phi_n: \Conf_{n+1}(M)\to \Conf_n(M)$ given by ``forget the point labeled $i$'', where $1\leq i\leq n+1$.  The problem is, as we've seen above, the maps $\phi_n$ typically do not induce isomorphisms on homology.   What to do?

Representation stability  allowed Church to analyze this situation.  It provided him with a language so that he could prove (Theorem~\ref{theorem:church}) that the maps $\phi_n$ stabilize.  The power of this point of view can be seen by applying Church's theorem (Theorem~\ref{theorem:church}) to the trivial representation $V(0)$, which gives (for $\dim(M)\geq 3$) that $d_{i,n}(0)$ is constant for $n\geq 2i$.  Transfer implies (see \eqref{eq:transfer1} above) that 
 $\dim_\C H_i(\UConf_{n}(M);\C)$ is constant for $n\geq 2i$, giving classical stability without maps between the spaces!  Church actually obtains the better stable range of $n>i$ by a more careful analysis.    
 
 One might notice that Church only obtains homological stability over $\Q$, while McDuff and Segal's theorem works over $\Z$.   One crucial place where $\Q$ is needed is the use of transfer.  However, this is not just an artifact of Church's proof, it is a feature of the situation: classical homological stability for 
 $\UConf_n(M)$ with $M$ closed is {\em false} for general closed manifolds $M$! For example,  $H_1(\UConf_n(S^2);\Z)=\Z/(2n-2)\Z$.    After Church's paper appeared, other proofs of homological stability over for $H_i(\UConf_n(M);\Q)$ were given by Randal-Williams \cite{RW} and then by Bendersky-Miller \cite{BM}.  
 
 By plugging other representations into Theorem~\ref{theorem:church}, Church deduces classical homological stability for a number of other colored configuration spaces.   The above discussion illustrates how representation stability can be used as a useful method to discover and prove classical homological stability theorems.   
 
 \subsection{Murnaghan's Theorem}
The stabilization of names of natural sequences of representations is not new.  The notation $V(\lambda)$ that we gave in \S\ref{section:configs} above goes back at least to the 1938 paper \cite{Mu} of Murnaghan, where he discovered the following theorem, first proved by Littlewood \cite{Li}  in 1957.

\begin{theorem}[{\bf Murnaghan's Theorem}] 
\label{theorem:Murnaghan}
For any two partitions $\lambda, \mu$ there is a finite set $P$ of partitions so that for all sufficiently large $n$:
\begin{equation}
\label{eq:Murnaghan}
V(\lambda)_n \otimes V(\mu)_n = \bigoplus_{\nu \in P} d_{\lambda,\mu}(\nu) V(\nu)_n
\end{equation}
for some non-negative integers $d_{\lambda\mu}(\nu)$.
\end{theorem}

The integers $d_{\lambda\mu}(\nu)$ are called {\em Kronecker coefficients}.  In his original paper \cite{Mu} Murnaghan computes the $d_{\lambda\mu}(\nu)$ explicitly for 58 of the simplest 
pairs $\mu,\nu$.    The study of Kronecker coefficients remains an active direction for research.  It is central to combinatorial representation theory and geometric complexity theory, among other areas.  
See, for example, \cite{BDO} and the references contained therein.  

One can deduce from Murnaghan's Theorem that the sequence $V(\lambda)_n \otimes V(\mu)_n$ is multiplicity stable in the sense of Definition~\ref{definition:repstab}; see \cite{CF1}.    In the following section we will describe a theory where Murnaghan's Theorem pops out as a structural feature of the theory.

\section{FI-modules}
\label{section:FI}

Representation stability for symmetric groups $S_n$ grew in power and applicability in \cite{CEF1}, where Thomas Church, Jordan Ellenberg and I developed a 
theory of  FI-modules.    

When looking at the sequence $H^i(\Conf_n(M);\C)$, we broke symmetry by only considering the map $\Conf_{n+1}(M)\to\Conf_n(M)$ given by ``forget the $(n+1)^{\rm st}$ point''. 
Of course there are really $n+1$ equally natural maps, given by ``forget the $j^{\rm th}$ point'' for $1\leq j\leq n+1$.   Taking cohomology switches the direction of 
arrows, and we have $n+1$ homomorphisms $H^i(\Conf_n(M);\C)\to H^i(\Conf_{n+1}(M);\C)$, 
each one corresponding to an injective map $\{1,\ldots ,n\}\to \{1,\ldots ,n+1\}$; namely, the injective map whose image misses $j$.    It is useful to consider all of these maps at once. This is the starting point for the study of FI-modules.

\subsection{FI-module basics} 
 An {\em FI-module} $V$ is a functor from the category FI of finite sets and injections to the category of modules over a fixed Noetherian ring $k$.   Thus to each set ${\bf n}:=\{1,\ldots ,n\}$ with $n$ elements the functor $V$ associates a $k$-module $V_n:=V(\bf n)$, and to each injection ${\bf m}\to {\bf n}$ the functor $V$ associates a linear map $V_m\to V_n$.  The set of self-injections ${\bf n}\to {\bf n}$ is the symmetric group $S_n$.    Thus an FI-module gives a sequence of $S_n$-representations $V_n$ and linear maps between them, one for each injection of finite sets:

\[\begin{array}{llllll}
\{1\}\longrightarrow&\{1,2\}\longrightarrow&\{1,2,3\}\longrightarrow&\cdots \longrightarrow& \{1,\ldots ,n\}\longrightarrow& \cdots \\
\ \circlearrowleft &\ \  \circlearrowleft &\ \ \ \  \circlearrowleft& & \ \ \ \ \ \circlearrowleft\\
\ S_1&\ \ S_2&\ \ \ \ S_3& &\ \ \ \ \ S_n \\
& & & & \\
& & & & \\
& & \text{\huge $V \downarrow$} & & \\
& & & & \\
& & & & \\
V_1\longrightarrow&V_2\longrightarrow&V_3\longrightarrow&\cdots \longrightarrow& V_n\longrightarrow& \cdots \\
\ \circlearrowleft &\  \circlearrowleft &\ \circlearrowleft& & \ \circlearrowleft\\
\ S_1&\ S_2&\ S_3& &\ S_n \\

\end{array}
\]

Of course each single horizontal arrow really represents many arrows, one for each injection between the corresponding finite sets.   Using functors from the category FI to study sequences of objects is not new: FI-spaces were known long ago (under different names) to homotopy theorists.  

A crucial observation is that one should think of an FI-module as a module in the classical sense.    
Many of the familiar notions from the theory of modules, such as submodule and quotient module, carry over to FI-modules in the obvious way: one performs the operations pointwise.  So, for example, $W$ is an FI-submodule of $V$ if $W_n\subset V_n$ for each $n\geq 1$.   One theme of \cite{CEF1} is that there is conceptual power in the encoding of this large amount of (potentially complicated) data into a single object $V$. 

The property of being an FI-module itself does not guarantee much structure.  One of the main insights in 
\cite{CEF1} was to find a finite generation condition that has strong implications but that one can also prove to hold in many examples. 

\begin{definition}[{\bf Finite generation}]
An FI-module $V$ is \emph{finitely generated} if there is a finite set $S$ of elements in $\coprod_iV_i$ so that no proper sub-FI-module of $V$ contains $S$.
\end{definition}

\begin{example}
Let $k[x_1,\ldots ,x_n]_{(3)}$ denote the vector space of homogeneous polynomials of degree $3$ in $n$ variables over a field $k$.  It is not hard to check that $\{1,\ldots ,n\}\mapsto k[x_1,\ldots ,x_n]_{(3)}$ is an FI-module.  We claim that  this FI-module is finitely generated by the elements 
$x_1^3, x_1^2x_2$ and $x_1x_2x_3$:

\[\begin{array}{lllll}
k[x_1]_{(3)}\longrightarrow&k[x_1,x_2]_{(3)}\longrightarrow&k[x_1,x_2,x_3]_{(3)}
\longrightarrow&k[x_1,x_2,x_3,x_4]_{(3)}
\longrightarrow\cdots &\\
& & & & \\
\ \fbox{$x_1^3$} & x_1^3,x_2^3&x_1^3,x_2^3,x_3^3&x_1^3,x_2^3,x_3^3,x_4^3 \\
&  \fbox{$ x_1^2x_2$}, x_2^2x_1  & x_1^2x_2, x_1^2x_3,x_2^2x_1&\ \ \ \  \ \ \ \vdots & \\
& & x_2^2x_3, x_3^2x_1,x_3^2x_2 &\\
&&\fbox{$x_1x_2x_3$}&

\end{array}
\]

\bigskip
Here we have written below each $k[x_1,\ldots ,x_n]_{(3)}$ its basis as a vector space.  Finite generation of the FI-module $k[x_1,\ldots ,x_n]_{(3)}$ is simply the fact that every vector in every 
$k[x_1,\ldots ,x_n]_{(3)}$ lies in the $k$-span of the set of vectors that can be obtained from the three boxed vectors by performing all possible morphisms, i.e.\ by changing the labels of the $x_i$.  In other words, there are, up to labeling and taking linear combinations, only three homogeneous degree three polynomials in any number $n\geq 3$ of variables: $x_1^3, x_1^2x_2$ and $x_1x_2x_3$.  Note that we need $n\geq 3$ to obtain all of the generators.  Similarly, $k[x_1,\ldots ,x_n]_{(87)}$ is finitely generated, but the full generating set appears only for $n\geq 87$.

\end{example}

The connection of FI-modules with representation stability is the following, proved in \cite{CEF1}.

\begin{theorem}[{\bf Finite generation vs.\ representation stable}]
\label{theorem:firepstab}
Let $V$ be an FI-module over a field $k$ of characteristic $0$.  Then $V$ is finitely generated if and only if $\{V_n\}$ is a representation stable sequence of $S_n$-representations with $\dim_kV_n<\infty$ for all $n$.
\end{theorem}

Theorem~\ref{theorem:firepstab} thus converts a somewhat complicated property about a sequence $V_n$ of representations into a single property -- finite generation -- of a single object $V$.    One example of the power of this viewpoint is the following.

\medskip
\noindent
{\bf Proof of Murnaghan's Theorem (Theorem~\ref{theorem:Murnaghan}): } Since $V(\lambda)$ and $V(\mu)$ are finitely generated FI-modules \cite[\S 2.8]{CEF1}, so is $V(\lambda)\otimes V(\mu)$.   Theorem~\ref{theorem:firepstab} implies that $V(\lambda)_n\otimes V(\mu)_n$ is representation stable, and so the theorem follows. \qed

\medskip
Thus a combinatorial theorem about an infinite list of numbers falls out of a basic structural property of FI-modules.

\subsection{Character polynomials} 
\label{section:charpoly}
One of the main discoveries of \cite{CEF1} is that character polynomials, studied by Frobenius but not so widely known today, are ubiquitous, and are an incredibly concise way to encode stability phenomena for sequences of $S_n$-representations. 

Fix the ground field $\C$.  Recall that the {\em character} of a representation $\rho:G\to\GL(V)$ of a finite group $G$ over $\C$ is defined 
to be the function $\chi_V:G\to \C$ given by \[\chi_V(g):=\text{Trace}(\rho(g)).\]   We view $\chi_V$ as an element of 
the vector space ${\cal C}(G)$ of {\em class functions} on $G$; that is, those functions that are constant on each conjugacy class in $G$.  A fundamental theorem in the representation theory of finite groups is that any $G$-representation is determined by its character: 

\begin{center}
{\it $\chi_V=\chi_W$ in ${\cal C}(G)$ if and only if $V\approx W$ as $G$-representations.}
\end{center}

For each $i\geq 1$ let $X_i\colon \coprod_nS_n\to \N$ be the class function defined by 
\[X_i(\sigma)= \mbox{number of $i$-cycles in the cycle decomposition of $\sigma$}.\]
A {\em character polynomial} is any polynomial $P\in\Q[X_1,X_2,\ldots ]$.   Such a polynomial gives a class function on all the $S_n$ at once.  The study of character polynomials goes back to work of Frobenius, Murnaghan, Specht,  and Macdonald; see, e.g.\ \cite[Example I.7.14]{Mac}).

It is easy to see for any fixed $n\geq 1$ that ${\cal C}(S_n)$ is spanned by character polynomials, so the character of any representation can be described by such a polynomial. 
For example, if $\C^n$ is the standard permutation representation of $S_n$ then the character $\chi_{\C^n}(\sigma)$ is the number of fixed points of $\sigma$, so $\chi_{\C^n}=X_1$ for any $n\geq 1$.  As another example, consider the $S_n$-representation $\bigwedge^2\C^n$.  Since $\sigma\cdot (e_i\wedge e_j)=\pm e_i\wedge e_j$ according to whether $\sigma$ contains $(i)(j)$ or $(ij)$, respectively, 
it follows that 

\[\chi_{\bigwedge^2\C^n}=\displaystyle{\binom{X_1}{2}-X_2}=\frac{1}{2}X_1^2-\frac{1}{2}X_1-X_2\]
for any $n\geq 1$.  These descriptions of characters are uniform in $n$.   On the other hand, if one fixes $r$ then for $n\gg r$ it is incredibly rare for an $S_n$-representation to be given by a character polynomial $P(X_1,\ldots ,X_r)$ depending only on cycles of length at most $r$.  A simple example is the sign representation: for $n\gg r$ one cannot determine the sign of an arbitrary $\sigma\in S_n$ just by looking at  cycles in $\sigma$ of length at most $r$. 

One of the main discoveries in \cite{CEF1} is that finitely-generated FI-modules in characteristic~0 
admit such a uniform description.

\begin{theorem}[{\bf Polynomiality of characters}] 
\label{thm:intro:persinomial}
Let $V$ be a finitely-generated FI-module over a field $k$ of characteristic~$0$.  Then the sequence of characters $\chi_{V_n}$ of the $S_n$-representations $V_n$ is {\em 
eventually polynomial}: there exists $N\geq 0$ and a polynomial $P(X_1,\ldots ,X_r)$ for some $r>0$ so that 
\[\chi_{V_n}=P(X_1,\ldots ,X_r) \ \ \ \mbox{for all $n\geq N$}\]
In particular $\dim_k(V_n)$ is a polynomial in $n$ for $n\geq N$.
\end{theorem}

The claim on $\dim_k(V_n)$ is obtained by noting that $$\dim_k(V_n)=\chi_{V_n}({\rm Id})=P(n,0,\ldots ,0).$$   
The fact that $\dim_k(V_n)$ is eventually a polynomial was was extended to the case ${\rm char}(k)>0$ in \cite{CEFN}.  In situations of interest one can often give explicit bounds on $r$ and $N$.  This converts the problem of finding all the characters $\chi_{V_n}$ into a concrete finite computation.   In some cases one can even get $N=0$.  
 
We again emphasize that the impact of Theorem~\ref{thm:intro:persinomial} comes not just from the fact that a single polynomial gives all characters of all $V_n$ with $n\gg 1$ at the same time, but it gives an extremely strong constraint on each individual $V_n$ for $n\gg r$, since $\chi_{V_n}=P(X_1,\ldots ,X_r)$ depends only on cycles of length at most $r$. 

\subsection{Examples/Applications} 
\label{section:fiexamples}
Part of the usefulness of finitely generated FI-modules is that they are common.   
This is illustrated in Table~\ref{table:examples}.   We define only a few of these examples here;  see  \cite{CEF1} for a detailed discussion.  

\begin{table}[h]
\centering
\begin{tabular}{ll}
\underline{FI-module $V=\{V_n\}$} \qquad\quad& \underline{Description} \\

&\\
$1. \ H^i(\Conf_n(M);\Q)$ & $\Conf_n(M)=$ configuration space of $n$ distinct ordered points \\
& on a connected, oriented manifold $M$, $\dim(M)>1$\\
&\\
2.\ $R^{(r)}_J(n)$
& $J=(j_1,\ldots j_r)$, $R^{(r)}(n)=\bigoplus_{J}R^{(r)}_J(n)$= $r$-diagonal \\
&   coinvariant 
algebra on $r$ sets of  $n$ variables \\

&\\
3.\ $H^i({\cal M}_{g,n};\Q)$ & ${\cal M}_{g,n}=$ moduli space of $n$-pointed genus $g\geq 2$ curves\ \ \\
&\\
4.\ $\cR^i({\cal M}_{g,n})$ & $i^{\text{th}}$ graded piece of the tautological ring of ${\cal M}_{g,n}$\   \\
&\\
5. ${\cal O}(X_{P,r}(n))_i$ & space of degree $i$ polynomials on the rank variety\\ 
& $X_{P,r}(n)$ of 
$n\times n$  matrices of $P$-rank $\leq r$ \\
&\\
6.\ $G(A_n/\Q)_i$ & degree $i$ part of the Bhargava--Satriano Galois closure\\ & of $A_n=\Q[x_1,\ldots,x_n]/(x_1,\ldots,x_n)^2$ \\
&\\
7.  $\langle H^1(\I_n;\Q)\rangle_{(i)}$ & degree $i$ part of the subalgebra of $H^\ast(\I_n;\Q)$ generated by \\
& $H^1(\I_n;\Q)$, where $\I_n=$ genus $n$ Torelli group\\
&\\
8. $H^i({\rm BDiff}_n(M);\Q)$ & ${\rm BDiff}_n(M)=$\ Classifying space of diffeos leaving a 
given set \\
&  of $n$ points invariant, for many manifolds $M$ (see \cite{JR2})\\
&\\
9.\ $\gr(\Gamma_n)_i$ &  $i$th graded piece of associated graded Lie algebra of many \\
& groups $\Gamma_n$, including $\I_n$, $\IA_n$ and the pure braid group $P_n$\\
\end{tabular}
\caption{Some examples of finitely generated FI-modules. Any parameter not equal to $n$ should be considered fixed and nonnegative. }
\label{table:examples}
\end{table}

\begin{theorem}[{\bf Finite generation}]
\label{theorem:examples:fg}
Each of the FI-modules (1)-(9) given in Table~\ref{table:examples} is finitely generated.
\end{theorem}

Items 3 and 8 of Theorem~\ref{theorem:examples:fg} are due to Jimenez Rolland \cite{JR1,JR2}; the other items are due to Church-Farb-Ellenberg \cite{CEF1}.

That each of (1)-(9) in Table~\ref{table:examples} is an FI-module is not difficult to prove.  More substantial is proving finite generation.  To do this one of course needs detailed information about the specific example.  In some of the cases this involves significant (but known) results; see below.  

Except for a few special (e.g.\ $M=\R^d$) and low-complexity  (i.e.\ small $i$, $d$, $g$, $J$, etc.) cases, explicit formulas for the characters (or even the dimensions) of the vector spaces  (1)-(9) of Table~\ref{table:examples} do not seem to be known, or even conjectured.  Exact computations may be quite difficult.  Applying Theorem~\ref{thm:intro:persinomial} and Theorem~\ref{theorem:examples:fg} to these examples gives us an answer, albeit a non-explicit one, in all cases.

\begin{theorem}[{\bf  Ubiquity of character polynomials}]
\label{thm:ubiquity}
For each of the sequences $V_n$ in Table~\ref{table:examples} there are numbers $N\geq 0, r\geq 1$ and a polynomial $P(X_1,\ldots ,X_r)$ so that 
\[\chi_{V_n}=P(X_1,\ldots ,X_r) \ \ \ \mbox{for all $n\geq N$}\]
In particular $\dim(V_n)$ is a polynomial in $n$ for $n\geq N$. 
\end{theorem} 

We emphasize that we are claiming eventual {\em equality} to a polynomial, not just polynomial growth.  As a contrasting example, if $\overline{\cal M}_{g,n}$ is the Deligne-Mumford compactification of the moduli space of $n$-pointed genus $g$ curves, the dimension of $H^2(\overline{\cal M}_{g,n};\Q)$ grows exponentially with $n$;  in particular the character of $H^2(\overline{\cal M}_{g,n};\Q)$ is not given by a character polynomial.  Although $V_n:=H^2(\overline{\cal M}_{g,n};\Q)$ is an FI-module, this FI-module is not finitely generated.

As an explicit example of Theorem~\ref{thm:ubiquity}, the character of the $S_n$-representation $H^2(\Conf_n(\C);\Q)$ is given for all $n\geq 0$  by the character polynomial
\begin{equation}
\label{eq:h2conf2}
\chi_{\displaystyle H^2(\Conf_n(\R^2);\C)}=2\binom{X_1}{3}+3\binom{X_1}{4}+\binom{X_1}{2}X_2-\binom{X_2}{2}-X_3-X_4.
\end{equation}

Note that for general finitely-generated FI-modules we only know such information for $n\gg 1$.   Recall from \eqref{eq:stabilizationinaction}  of  \S\ref{section:repstab} how decomposition into irreducibles of $H^2(\Conf_n(\R^2);\C)$ was shown to change with $n$, only to stabilize once $n\geq 7$.  I encourage the reader to try to see this via \eqref{eq:h2conf2}, which 
holds for all $n\geq 0$.  

Although we can sometimes give explicit upper bounds on their degree, the polynomials produced by 
Theorem~\ref{thm:ubiquity} are known explicitly in only a few special cases.   Thus the following is one of the main open problems in this direction.

\begin{problem}
Compute the polynomials $P(X_1,\ldots ,X_r)$ produced by Theorem~\ref{thm:ubiquity}.
\end{problem}

One difficulty in solving this problem is that, in many examples, the proof of finite generation of the corresponding FI-module uses a Noetherian property (see below), and the proof of this property is not effective.

\subsection{The Noetherian property} The following theorem,  joint work with Thomas Church, Jordan Ellenberg, and Rohit Nagpal, is central to the theory of FI-modules; it is perhaps the most useful general tool for proving that a given FI-module is finitely generated.  

\begin{theorem}[{\bf Noetherian property}]
\label{thm:noetherian}
Let $V$ be a finitely-generated FI-module over a Noetherian ring $k$.  Then any sub-FI-module of $V$ is finitely generated.
\end{theorem}

Theorem~\ref{thm:noetherian} was proved in this generality by Church-Ellenberg-Farb-Nagpal \cite{CEFN}.   For fields $k$ of characteristic $0$ it was proved earlier by Church-Ellenberg-Farb \cite[Theorem~2.60]{CEF1} and by Snowden \cite[Theorem~2.3]{Sn}, who actually proved a version (in a different language) for modules for many twisted commutative algebras, of which FI-modules are an example.   The Noetherian property for FI-modules over fields $k$ of positive characteristic is crucial 
for the study of the cohomology of congruence subgroups from this point of view; see \S\ref{section:congruence} below.  L\"{u}ck proved a version of Theorem~\ref{thm:noetherian} for finite categories in \cite{Lu}, but since FI is infinite these do not occur in our context.

One can see how Theorem~\ref{thm:noetherian} is used in practice via the following.

\begin{theorem}
\label{theorem:spectral:sequence}
Suppose $E_*^{p,q}$ is a first-quadrant spectral sequence of FI-modules over a Noetherian ring $k$, and that $E_*^{p,q}$ converges to an FI-module $H^{p+q}(X;k)$.  If the FI-module $E_2^{p,q}$ is finitely generated for each $p,q\geq 0$, then the FI-module $H^i(X;k)$ is finitely generated for each fixed $i\geq 0$.
\end{theorem}

See \cite{Ch} and \cite{CEF1} for earlier versions of Theorem~\ref{theorem:spectral:sequence}, and 
the paper \cite{JR2}, where Jimenez Rolland gives explicit bounds on the stability degree, etc.

Spectral sequences as in Theorem~\ref{theorem:spectral:sequence} arise in many computations.  For example,  following Cohen-Taylor \cite{CT} and Totaro \cite{To}, one computes $H^i(\Conf_n(M);k)$ by using the Leray spectral sequence for the natural inclusion $\Conf_n(M)\to M^n$.   As $n$ varies we obtain a sequence of spectral sequences, one for each $n$.  In fact this gives a spectral sequence of FI-modules.    Another example is the computation of the homology of congruence subgroups (see \S\ref{section:congruence} below).

The proof of Theorem~\ref{theorem:spectral:sequence} is that, while we have no idea what the differentials might be, or at which page the spectral sequence stabilizes (and this may depend on $n$),  the terms $E^{p,q}_\infty$ are obtained from the $E_2^{p,q}$ terms by repeatedly taking submodules and quotient modules.  Since the property of finite generation for an FI-module is preserved by taking submodules (by the Noetherian property Theorem~\ref{thm:noetherian}) and quotients, then if $E_2^{p,q}$ is finitely generated so is $E_j^{p,q}$ for every $j\geq 2$. 

\subsection{Some remarks on the general theory}  There are many other aspects of the general theory of FI-modules that I am not describing due to lack of space.  This includes a more quantitative version of the theory, with notions such as stability degree and weight of an FI-module, allowing for 
explicit estimates on stable ranges and degrees of character polynomials.  Co-FI-modules are useful 
when one has maps going the wrong way.    Also useful are FI-spaces, FI-varieties, and FI-hyperplane arrangements (see \cite{CEF2} for the latter); these give FI-modules by applying the (co)homology functor.  Church-Putman \cite{CP} have developed the theory of FI-groups in order to prove a kind of relative finite generation theorem  in group theory; they apply this in \cite{CP} to certain subgroups of Torelli groups.  Sam and Snowden have given in \cite{SS} a more detailed analysis of the algebraic structure of the category of FI-modules in characteristic~0.

\section{Combinatorial statistics for varieties over finite fields}
\label{section:varieties}

In \cite{CEF2} we exposed a close connection between representation stability in cohomology and the stability of various combinatorial statistics for polynomials over finite fields and for maximal tori in $\GL_n(\F_q)$. We now give a brief sketch of how this works.

\subsection{The space of polynomials over $\F_q$} Consider the following basic questions: how many square-free (i.e.\ having no repeated roots), degree $n$ monic polynomials in $\F_q[T]$ are there? How many linear factors does one expect such a polynomial to have? factors of degree $d$? What is the variance of this expectation?   

If one fixes $q$ and allows $n$ to increase, something interesting happens.  A good example of what I'd like to describe is the expected {\em quadratic excess} of a polynomial in $\F_q[T]$; that is, the expected difference of the number of reducible quadratic factors and the number of irreducible quadratic factors.   This number can be computed by adding up the quadratic excess of each degree $n$, monic square-free polynomial in $\F_q[T]$  and then dividing by the total number 
$q^n-q^{n-1}$ of such polynomials.  Here are some values for small $n$:

 \begin{align*}
&\text{total: }&&\text{expectation: }\\
n=3:\qquad&q^2-q&&\frac{1}{q}\\
n=4:\qquad&q^3-3q^2+2q&&\frac{1}{q}-\frac{2}{q^2}\\
n=5:\qquad&q^4-4q^3+5q^2-2q&&\frac{1}{q}-\frac{3}{q^2}+\frac{2}{q^3}\\
n=6:\qquad&q^5-4q^4+7q^3-7q^2+3q&&\frac{1}{q}-\frac{3}{q^2}+\frac{4}{q^3}-\frac{3}{q^4}\\
n=7:\qquad&q^6-4q^5+7q^4-8q^3+8q^2-4q&&\frac{1}{q}-\frac{3}{q^2}+\frac{4}{q^3}-\frac{4}{q^4}+\frac{4}{q^5}\\
n=8:\qquad&q^7-4q^6+7q^5-8q^4+9q^3-10q^2+4q\quad&&\frac{1}{q}-\frac{3}{q^2}+\frac{4}{q^3}-\frac{4}{q^4}+\frac{5}{q^5}-\frac{5}{q^6}
\end{align*}

Notice that in each column in the counts above, the coefficient changes as $n$ increases, 
until $n$ is sufficiently large, and then this coefficient stabilizes.  For example the third column gives coefficients $0,2,5,7,7,7,\ldots $   Theorem~\ref{theorem:polystat} below implies that these formulas must converge term-by-term to a limit.  A somewhat involved computation (see below) allowed us in \cite{CEF2} to compute this limit as:
\begin{equation}
\label{eq:intro:limitqe}
q^{n-1}-4q^{n-2}+7q^{n-3}-8q^{n-4}+\cdots\qquad\quad\text{and}\quad\qquad\frac{1}{q}-\frac{3}{q^2}+\frac{4}{q^3}-\frac{4}{q^4}+\cdots,
\end{equation}

This numerical stabilization is a reflection of something deeper.  To explain, consider the space $Z_n$ of all monic, square-free, degree $n$ polynomials with coefficients in the  finite field $\F_q$.    Recall that the {\em discriminant} $\Delta_n\in\Z[x_1,\ldots ,x_{n-1}]$ is a polynomial with the  property that an arbitrary monic polynomial $f(z)=z^n+a_{n-1}z^{n-1}+\cdots +a_1z_1+a_0\in\C[z]$ is square-free if and only if $\Delta_n(a_0,\ldots a_{n-1})\neq 0$.   Thus $Z_n$ is a complex algebraic variety. For example 
\[
\begin{array}{c}
Z_2=\{z^2+bz+c\in\C[z] : b^2-4c\neq 0\}\\
\\
\text{and}\\
\\
Z_3=\{z^3+bz^2+cz+d\in \C[z]: b^2c^2-4c^3-4b^3d-27d^2+18bcd\neq 0\}.
\end{array}\]

The complex variety $Z_n$ is also an algebraic variety over the finite field $\F_q$ for any prime power $q$.  The set of $\F_q$-points $Z_n(\F_q)$ is exactly the set of monic, square-free, degree $n$ polynomials in $\F_q[T]$.    From this point of view we should think of the original complex algebraic variety as the complex points $Z_n(\C)$. There is a remarkable 
relationship between $Z_n(\C)$ and $Z_n(\F_q)$, given by the Grothendieck--Lefschetz fixed point theorem in \'{e}tale cohomology, which we now explain.

\subsection{The Grothendieck-Lefschetz formula}
It is a fundamental observation of Weil that for an algebraic variety $Z$ defined over $\F_q$, one can realize $Z(\F_q)$ as the fixed points of a dynamical system, as follows.  Denote by $\Fqbar$ the algebraic closure of $\F_q$.  The {\em geometric Frobenius} morphism 
$\Frob_q\colon Z(\Fqbar)\to Z(\Fqbar)$ acts (in an affine chart) on the coordinates of $Z$ by 
$x\mapsto x^q$.  Fermat's Little Theorem implies that 
\[Z(\F_q)=\Fix [(\Frob_q: Z(\Fqbar)\to Z(\Fqbar)].\]
In the case of the varieties $Z_n$ that we are considering, the Grothendieck-Lefschetz fixed point theorem takes the form: 
\begin{equation}
\label{eq:groth1}
|Z_n(\F_q)|=\sum_{f\in \Fix(\Frob_q)}\!\! \! \! \! \! 1=\sum_i (-1)^iq^{n-i}\dim_\C H^i(Z_n(\C);\C)=q^n-q^{n-1}
\end{equation}
where the last equality comes from the theorem of Arnol'd that $H^i(Z_n(\C);\C)=0$ unless $i=0,1$, in which case it is $\C$.    

We want to compute more subtle counts than just $|Z_n(\F_q)|$.  To this end, we can {\em weight} the points of $Z_n(\F_q)=\Fix(\Frob_q)$, as follows.  For each $f\in \Fix(\Frob_q)=Z_n(\F_q)$ the map $\Frob_q$ permutes the set 
\[\roots (f):=\{y\in\Fqbar: f(y)=0\}\] giving a conjugacy class $\sigma_f$ in $S_n$.  Thus 
$X_i(\sigma_f)$ is well defined.   Let $d_i(f)$ denote the number of irreducible (over $\F_q$) degree $i$ factors of $f$.  A crucial observation is that for any $i\geq 1$:
\begin{equation}
\label{eq:keyequality1}
X_i(\sigma_f)= d_i(f).
\end{equation}

Any $P\in\C[x_1,\ldots ,x_r]$ (here $r\geq 1$ is arbitrary) determines a character polynomial $P(X_1,\ldots, X_r)$ (cf.\ \S\ref{section:charpoly} above).  The polynomial $P$ gives a way to weight points $f\in Z_n(\F_q)$  via 
\[P(f):= P(d_1(f),\ldots ,d_r(f))=P(X_1(\sigma_f),\ldots ,X_r(\sigma_f)).\]
So for example $P(f)=d_1(f)$ counts the number of linear factors of $f$ (i.e.\ the number of roots of $f$ that lie in $\F_q$), and $P(f)=d_2(f)-d_1(f)^2$ counts the quadratic excess of $f$.  In general we call such a $P$ a {\em polynomial statistic}.  The expected value of $P(f)$ for $f\in Z_n(\F_q)$ is then given by $[\sum_{f\in Z_n(\F_q)}P(f)]/(q^n-q^{n-1})$.

\para{\boldmath$H^i(\Conf_n(\C);\C)$ enters the picture} 
Computing this expectation for a given $P$ is where $H^i(\Conf_n(\C);\C)$ comes in.  We can identify $Z_n(\C)$ with the space $\UConf_n(\C)=\Conf_n(\C)/S_n$ of unordered $n$-tuples of distinct points in $\C$ via the bijection that sends $f\in Z_n(\C)$ to its set of roots.    We thus have a covering $\Conf_n(\C)\to Z_n(\C)$ of algebraic varieties, with deck group $S_n$.  

Now, it's something of a long story, and there are a number of technical details to worry about, but the theory of  \'{e}tale cohomology and the twisted Grothendieck-Lefschetz formula, together with work of Lehrer \cite{Leh}, who proved that this machinery can be applied in this case,  can be used to give the following theorem of \cite{CEF2}.  Let 
$\langle \phi,\psi\rangle_{S_n}:=\sum_{\sigma\in S_n}\phi(\sigma)\overline{\psi(\sigma)}$ be the standard inner product on the space of $\C$-valued functions on $S_n$.

\begin{theorem}[{\bf Twisted Grothendieck--Lefschetz for \boldmath$Z_n$}]
For each prime power $q$, each positive integer $n$, and each character polynomial $P$, we have 
\begin{equation}
\sum_{\displaystyle f \in Z_n(\F_q)} \!\!\!\!\!\!\! P(f) = \sum_{i=0}^n (-1)^i q^{n-i} \big\langle P,  \chi_{\displaystyle H^i(\Conf_n(\C);\C)}\big\rangle_{S_n}.
\end{equation}
\label{pr:exactcount}
\end{theorem}

For example, when $P = 1$ the inner product $\langle P, H^i(\Conf_n(\C);\C) \rangle$ is the  multiplicity of the trivial $S_n$-representation in $H^i(\Conf_n(\C);\C)$, which by transfer is the  dimension of $H^i(Z_n(\C);\C)$, giving the formula \eqref{eq:groth1} above.   Theorem~\ref{pr:exactcount} tells us that we can compute various weighted point counts on $Z_n(\F_q)$ if we understand the cohomology of the $S_n$-cover $\Conf_n(\C)$ of $Z_n(\C)$ as an $S_n$-representation. 

\subsection{Representation stability and Grothendieck-Lefschetz}
\label{section:count3}
 Now we can bring representation stability into the picture.  According to Theorem~\ref{thm:ubiquity}, for each $i\geq 0$ the character of $H^i(\Conf_n(\C);\C)$ is given by a character polynomial 
 for $n\gg 1$ (actually in this case it holds for all $n\geq 1$).  The inner product of two character polynomials is, for $n\gg 1$, constant. Keeping track of stable ranges, and defining $\deg P$ as usual but with $\deg x_k=k$, we deduce in \cite{CEFN} the following.

\begin{theorem}[{\bf Stability of polynomial statistics}]
\label{theorem:polystat}
For any polynomial $P\in\Q[x_1,x_2,\ldots ]$, the limit 
 \[\big\langle P,H^i(\Conf_{\bullet}(\C);\C)\big\rangle :=
 \lim_{n\to \infty}\big\langle P, \chi_{\displaystyle H^i(\Conf_n(\C);\C)}\big\rangle_{S_n}\] 
 exists; in fact, this sequence is constant for $n\geq 2i+\deg P$.  Furthermore, for each prime power $q$:
\begin{equation}
\label{eq:stable3}
\lim_{n \ra \infty} q^{-n} \sum_{f \in Z_n(\F_q)}  \!\!\!\!\!\! P(f) = \sum_{i=0}^\infty (-1)^i
\big\langle P,\chi_{\displaystyle H^i(\Conf_{\bullet}(\C);\C)}\big\rangle q^{-i}
\end{equation}
In particular, both the limit on the left and the series on the right in \eqref{eq:stable3} converge, and they converge to the same limit.
\end{theorem}

Plugging $P=\binom{X_1}{2}-X_2$ into Theorem~\ref{theorem:polystat} gives the stable formula 
\eqref{eq:intro:limitqe} for quadratic excess of a square-free degree $n$ polynomial in $\F_q[T]$.    The limiting values of other polynomial 
statistics $P$ are computed in \cite{CEF2}; some of these are given in Table~\ref{table:stats} below.  One can actually apply Equation~\eqref{eq:stable3} of Theorem~\ref{theorem:polystat} in reverse, using number theory to compute the left-hand side in order to determine the right-hand side, as we do in \S 4.3 of \cite{CEF2}.

The above method is applied in \cite{CEF2} in the same way to a different counting problem.  Consider the complex algebraic variety of ordered $n$-frames in $\C^n$: 
\[Z_n(\C)=\big\{(L_1,\ldots,L_n)\,\big|\,L_i \text{ a line in }\C^n,\ \ L_1,\ldots,L_n \text{ linearly independent}\big\}.\]
The group $S_n$ acts on $Z_n(\C)$ via $\sigma\cdot L_i=L_{\sigma(i)}$.  The quotient $Z_n(\C)/S_n$ is also an algebraic variety, and its $\F_q$-points parametrize the set of maximal tori in the finite group $\GL_n\F_q$.     In analogy with the case of square-free polynomials, each $P\in\C[X_1,\ldots ,X_r]$ counts maximal tori in $\GL_n\F_q$ with different weights.  Since $H^i(Z_n(\C);\C)$ is known to be representation stable (essentially by a theorem of Kraskiewicz-Weyman, Lustig,  and Stanley - see \S 7.1 of \cite{CF1}), we can apply an analogue of Theorem~\ref{theorem:polystat} in this context to compute this weighted point count.   

Table~\ref{table:stats} lists some examples of specific asymptotics that are computed in \cite{CEF2} using this method.    The formulas in each column are obtained from Theorem~\ref{theorem:polystat} (and its analogue for $Z_n(\C)$) with $P=1$, $P=X_1$, $P=\binom{X_1}{2}-X_2^2$, the character 
$\chi_{\text{sign}}$ of the sign representation, and the characteristic function $\chi_{n{\rm cyc}}$ of the $n$-cycle, respectively.  Note that the latter two are not character polynomials.

\begin{table}[h]
\centering
\begin{tabular}{lll}
&Counting theorem for\qquad\qquad\qquad\ &Counting theorem for\\
\underline{\ P\ }&\underline{squarefree polys in $\F_q[T]$}&\underline{maximal tori in $\GL_n \F_q$}\\[14pt]

$1$&\# of degree $n$ squarefree& \# of maximal tori in $\GL_n\F_q$\\
&polynomials $=q^n-q^{n-1}$&  (both split and non-split) $=q^{n^2-n}$\\[12pt]

$x_1$&expected \# of linear factors& expected \# of eigenvectors in $\F_q^n$\\
&$=1-\frac{1}{q}+\frac{1}{q^2}-\frac{1}{q^3}+\cdots \pm \frac{1}{q^{n-2}}$& $=1+\frac{1}{q}+\frac{1}{q^2}+\cdots +\frac{1}{q^{n-1}}$\\[12pt]

$\binom{x_1}{2}-x_2^2$&expected excess of \emph{reducible} & expected excess of \emph{reducible} \\
&vs. \emph{irreducible} quadratic factors &vs. \emph{irreducible} dim-2 subtori\\
&\ \ $\to\  \frac{1}{q}-\frac{3}{q^2}+\frac{4}{q^3}-\frac{4}{q^4}$&\ \ $\to\  \frac{1}{q}+\frac{1}{q^2}+\frac{2}{q^3}+\frac{2}{q^4}$\\
&\qquad$\ \ +\frac{5}{q^5}-\frac{7}{q^6}+\frac{8}{q^7}-\frac{8}{q^8}+\cdots$&\qquad$\ \ +\frac{3}{q^5}+\frac{3}{q^6}+\frac{4}{q^7}+\frac{4}{q^8}+\cdots$\\
&as $n\to \infty$&as $n\to \infty$\\[12pt]

$\chi_{\text{sign}}$&discriminant of random squarefree&\# of irreducible factors is more\\
&polynomial is equidistributed in $\F_q^\times$&likely to be $\equiv n\bmod{2}$ than not,\\
&between residues and nonresidues&with bias $\sqrt{\text{\# of tori}}$\\[12pt]

$\chi_{n{\rm cyc}}$&Prime Number Theorem for $\F_q[T]$:&\# of irreducible maximal tori\\
& \# of irreducible polynomials &$=\frac{q^{\binom{n}{2}}}{n}(q-1)(q^2-1)\cdots(q^{n-1}-1)$\\
&\quad$=\sum_{d|n}\frac{\mu(n/d)}{n}q^{d}\sim \frac{q^n}{n} $&\qquad\qquad$\sim c\cdot \frac{q^{n^2-n}}{n}$
\end{tabular}
\caption{Some asymptotics from \cite{CEF2}, computed using Theorem~\ref{theorem:polystat} and its $Z_n(\C)$ analogue.}
\label{table:stats}
\end{table}

The formulas for square-free polynomials in Table~\ref{table:stats} can be proved by direct means, for example using analytic number theory (e.g. weighted $L$-functions).  In contrast,  the formulas for maximal tori in $\GL_n \F_q$  may be known but are not so easy to prove. For example, the formula for the number of maximal tori in $\GL_n\F_q$ is a well-known theorem of Steinberg; proofs using the Grothendieck--Lefschetz formula have been given by  Lehrer and Srinivasan (see e.g.\ \cite{Sr}).  Regardless, a central message of \cite{CEF2}  is that representation stability provides a single underlying mechanism for all such formulas.

%

\begin{remark}[{\bf Stable range vs.\ rate of convergence}]
The dictionary between representation stability and stability of point-counts goes one level deeper.    
One can of course ask how the formulas in Table~\ref{table:stats} converge.  As discussed in \cite{CEF2}, the speed of convergence of any such formula depends on the stable range of the corresponding representation stability problem. For example, let $L$ denote the limit of each side of Equation~\eqref{eq:stable3}.  
The fact that $\langle \chi_P,H^i(\Conf_n(\C))\rangle_{S_n}$ is stable with stable range  $n\geq 2i+\deg P$  can be used to deduce that 
 \[q^{-n}\sum_{f(T)\in Z_n(\F_q)}  \!\!\!\!\!\!\!\!\!  P(f)=L+O(q^{(\deg P-n)/2})=L+O(q^{-n/2}).\] 
We thus have a {\em power-saving bound} on the error term.  See \cite{CEF2} for more details.\end{remark}

\section{FI-modules in characteristic $p$}
\label{section:congruence}

The Noetherian property for FI-modules was extended from fields of characteristic $0$ to arbitrary Noetherian rings by Church-Ellenberg-Farb-Nagpal \cite{CEFN}.  The proof is significantly more difficult in this case, and new ideas were needed.  Indeed, \cite{CEFN} brought in more categorical and homological methods into the theory, for example with a homological reformulation of finite generation, and a certain shift functor that plays a crucial role.  This line of ideas has culminated in the recent theory of FI-homology of Church-Ellenberg \cite{CE}, which is an exciting and powerful new tool.  

One reason that we care about characteristic $p>0$ is that in some examples this case contains most of the information. As an example, let $K$ be a number field with ring of integers ${\cal O}_K$, and let $\p\subset {\cal O}_K$ be any proper ideal.  
Define the {\em congruence subgroup} $\Gamma_n(\p)\subset \GL_n({\cal O}_k)$ to be 
\[\Gamma_n(\p):=\text{kernel}[\GL_n({\cal O}_k)\to \GL_n({\cal O}_k/\p)].\]

As shown by Charney \cite{Cha}, when one considers coefficients localized at $\p$ then $\Gamma_n(\p)$ and $\GL_n({\cal O}_K)$ have the same homology.  
Thus the interesting new information about $\Gamma_n(\p)$ comes via coefficients $k$ with ${\rm char}(k)=p>0$.  While $H_i(\Gamma_n(\p);k)$ is most naturally 
a representation of $\SL_n({\cal O}_k/\p)$, one can 
restrict this action to a copy of $S_n$, and show that $H_i(\Gamma_n(\p);k)$ is an FI-module.   
The Noetherian condition is crucial for proving that this FI-module is finitely generated, since the proof uses a spectral sequence argument (see below).  

The proof of Theorem~\ref{thm:intro:persinomial}, that for a finitely-generated FI-module $V$ the character $\chi_{V_n}$ is  a character polynomial for $n\gg 1$, 
works only over a field $k$ with ${\rm char}(k)=0$.   However, for fields $k$ 
with ${\rm char}(k)>0$, we were still able to prove \cite[Theorem B]{CEFN} that there is a polynomial $P\in\Q[T]$ so that $\dim_k(V_n)=P(n)$ for all $n\gg 1$.  Following the approach of Putman \cite{Pu}, we were able to apply this to 
$H_i(\Gamma_n(\p);k)$, giving the following theorem, first proved by Putman \cite{Pu} for fields of large characteristic.

\begin{theorem}[{\bf mod \boldmath$p$ Betti numbers of congruence subgroups}]
\label{thm:congruencepoly}
Let $K$ be a number field,  ${\cal O}_K$ its ring of integers, and $\p\subsetneq {\cal O}_K$ any proper ideal. For any $i\geq 0$ and any field $k$, there exists a polynomial $P(T)=P_{\p,i,k}(T)\in\Q[T]$ so that 
for all sufficiently large $n$, 
\[\dim_k H_i(\Gamma_n(\p);k) = P(n).\]
\end{theorem}

The exact numbers $\dim_k H_i(\Gamma_n(\p);k)$ for $i>1$ are known in very few cases, even for the simplest case $K=\Q, \p=(p), k=\F_p$.  Frank Calegari~\cite{Ca} has recently determined the rate of growth of the mod $p$ Betti numbers of the level $p^d$ congruence subgroup of $\SL_n({\cal O}_K)$. He proves for example in \cite[Lemma~3.5]{Ca} that for $p\geq 5,d\geq 1$: 
\[\dim_{\F_p}H_i(\Gamma_n(p^d);\F_p)=\binom{n^2-1}{i}+O(n^{2i-4}).\]
Calegari's result tells us the leading term of the polynomial guaranteed by Theorem~\ref{thm:congruencepoly}.  It should be noted that Calegari's proof uses (Putman's version of) 
Theorem~\ref{thm:congruencepoly}.

\begin{problem}[\cite{CEFN}]
Compute the polynomials $P_{\p,i,k}\in\Q[T]$ given by Theorem~\ref{thm:congruencepoly}.    Do the Brauer characters of $H_i(\Gamma_n(\p);k)$, or indeed of an arbitrary finitely-generated FI-module over a finite field $k$ with ${\rm char}(k)>0$, exhibit polynomial behavior in $n$ for $n\gg 1$?
\end{problem}

The more categorical setup in \cite{CEFN} allowed us to find an inductive description for any finitely generated FI-module.

\begin{theorem}[{\bf Inductive description of f.g. FI-modules}]
\label{thm:inductive}
Let $V$ be a finitely-generated FI-module over a Noetherian ring $R$. Then there exists some $N\geq 0$ such that for all $n\in \N$, there is an isomorphism of $S_n$-representations:
\begin{equation}
\label{eq:colim}
V_n\approx\ \varinjlim V(S)
\end{equation}
where the direct limit is taken over the poset of subsets $S\subset\{1,\ldots ,n\}$ with $|S|\leq N$.
\end{theorem}

The condition \eqref{eq:colim} in Theorem~\ref{thm:inductive} can be viewed as a reformulation of Putman's central stability condition \cite[\S1]{Pu}.   

Since we proved in \cite{CEFN} that $H_m(\Gamma_n(\p);k)$ is a finitely generated FI-module, Theorem~\ref{thm:inductive} thus gives the following inductive presentation of $H_m(\Gamma_n(\p);\Z)$.  Let $\Gamma_{n-1}^{(i)}(\p)$ with $1\leq i\leq n$ denote the $n$ standard subgroups  of $\Gamma_n(\p)$ isomorphic to $\Gamma_{n-1}(\p)$.  Let $\Gamma_{n-2}^{(i,j)}(\p):=\Gamma_{n-1}^{(i)}(\p)\cap \Gamma_{n-1}^{(j)}(\p)$.  As the notation suggests, each $\Gamma_{n-2}^{(i,j)}(\p)$ is isomorphic to $\Gamma_{n-2}(\p)$. As with the Mayer-Vietoris sequence, the difference of the two inclusions gives a map 
\[H_m(\Gamma^{(i,j)}_{n-2}(\p))\to H_m(\Gamma^{(i)}_{n-1}(\p))\oplus H_m(\Gamma^{(j)}_{n-1}(\p))\]  whose image vanishes in $H_m(\Gamma_n(\p))$.  A version of the following theorem for coefficients in a sufficiently large finite field was first proved by Putman \cite{Pu}. 

\begin{theorem}[{\bf A presentation for $H_m(\Gamma_n(\p);\Z)$}]
\label{thm:presentation}
Let $K$ be a number field, let ${\cal O}_K$ be its ring of integers, and let $\p$ be a proper ideal in ${\cal O}_K$. Fix $m\geq 0$. Then for all sufficiently large $n$, 
\[H_m(\Gamma_n(\p);\Z)\simeq\frac{\ \ \ \bigoplus_{i=1}^n H_m(\Gamma_{n-1}^{(i)}(\p);\Z)}{\im \bigoplus_{i< j} H_m(\Gamma_{n-2}^{(i,j)}(\p);\Z)}.\]
\end{theorem}

We think of Theorem~\ref{thm:presentation} as giving a presentation for $H_m(\Gamma_n(\p);\Z)$, with copies of $H_m(\Gamma_{n-1}(\p);\Z)$ as generators and copies of 
$H_m(\Gamma_{n-2}(\p);\Z)$ as relations.   Theorem~\ref{thm:inductive} is applied in \cite{CEFN} to give a similar description for $H_m(\Conf_n(M);\Z)$ and for graded pieces of diagonal coinvariant algebras.    Nagpal \cite{Na} has recently extended this point of view considerably, and has applied it to prove that the groups $H_m(\UConf_n(M);\F_p)$ are periodic in $n$.

\section{Representation stability for other sequences of representations}
\label{section:other:examples}

In this paper we focused our attention on sequences $V_n$ of $S_n$-representations.  This is just one of the examples from \cite{CF1}, where we introduced and studied representation stability (and variations) for other families $G_n$ 
of groups whose representation theory has a consistent naming
system.  Examples include $G_n=\GL_n\Q$, $\Sp_{2g}\Q$ and the hyperoctahedral
groups.  We also explored the case of modular representations of algebraic groups over finite fields, where instead of stability we found representation periodicity.   The reader is referred to \cite{CF1} for precise definitions and many examples.

I would like to illustrate here how these kinds of examples arise.  For brevity let's stick to the calculation of group homology.  Here is the general setup.  Let $\Gamma$ be a group with normal subgroup $N$ and quotient $A:=\Gamma/N$.  The conjugation action of $\Gamma$ on $N$
induces a $\Gamma$--action on $H_i(N;R)$ for any coefficient ring $R$. This action factors through an
$A$--action on $H_i(N,R)$, making $H_i(N,R)$ into an $A$--module.

The structure of $H_i(N,R)$ as an
$A$--module encodes fine information.  For example, the
transfer isomorphism shows that when $A$ is finite and $R=\Q$, the space
$H_i(\Gamma;\Q)$ appears precisely as the subspace of $A$--fixed
vectors in $H_i(N;\Q)$.  But there are typically many other summands,
and knowing the representation theory of $A$ (over $R$) gives us a
language with which to access these.

There are many natural examples of families $\Gamma_n$ of this type, with normal subgroups $N_n$ and quotients $A_n$.  Table~\ref{table:other} summarizes some examples that fit into this
framework.  

\begin{table}[h]
\label{table:other}
\centering
\begin{tabular}{c|c|c|c|c}
  kernel $N_n$ & group $\Gamma_n$ & acts on & quotient $A_n$
  & $H_1(N,R)$ for big $n$\\
  \hline 
  & & & & \\
pure braid group $P_n$ & braid group $B_n$ & $\{1,\ldots ,n\}$ & $S_n$ & ${\rm Sym}^2V_n/V_n$\\
  & & & & \\
  Torelli group $\I_n$ & mapping class&$H_1(\Sigma_n,\Z)$&$\Sp_{2n}\Z$& 
  $\bwedge^3V_n/V_n$ \\
  & group $\Mod_n$ & & & \\
  & & & & \\
  ${\rm IAut}(F_n)$&$\Aut(F_n)$&$H_1(F_n,\Z)$&$\GL_n\Z$
  &$V_n^\ast\otimes \bwedge^2V_n$\\
  & & & & \\
  congruence&$\SL_n\Z$&$\F_p^n$&$\SL_n\F_p$&$\fsl_n\F_p$\\
  subgroup $\Gamma_n(p)$& & & &\\
  & & & & \\
  level $p$ subgroup&$\Mod_n$&$H_1(\Sigma_n;\F_p)$&$\Sp_{2n}\F_p$&
  $\bwedge^3 V_n/V_n\oplus\ \fsp_{2n}\F_p$\\
  $\Mod_n(p)$&&&&
  \end{tabular}
  \caption{Some natural sequences of representations.}
  \end{table}

In each case the group $N_n$ arises as the kernel of a natural $\Gamma_n$-action.  
Each example is explained in detail in \cite{CF1}.    Here $R=\Q$ in the first three examples, $R=\F_p$ in the
fourth and fifth, and $V_n$ stands in each case for the ``standard
representation'' of $A_n$. In the last example $p$ is an odd prime. 

In each of the examples in Table~\ref{table:other}, the groups $\Gamma$ are known to
satisfy classical homological stability.  In contrast, the rightmost
column of Table~\ref{table:other} shows that none of the groups $N$ satisfies homological
stability, even in dimension 1.  In fact, except for the example of $P_n$,  
very little is known about the $A_n$--module $H_i(N_n,R)$ for
$i>1$, and indeed it is not clear if there is a nice closed form
description of these homology groups.  However, the appearance of some
kind of stability can already be seen in the rightmost column, as
the names of the irreducible composition factors of these $A_n$--modules are
constant for large enough $n$; this is discussed in detail in \cite{CF1}.

A crucial common property of the examples in Table~\ref{table:other} is that each of the sequences $A_n$ has an inherent stability in the naming of its irreducible algebraic
representations over $R$.  For example, an irreducible algebraic
representation of $\SL_n\Q$ is determined by its highest weight vector,
and these vectors may be described uniformly without reference to
$n$. For example, for $\SL_n\Q$ the irreducible representation
$V(L_1+L_2+L_3)$ with highest weight $L_1+L_2+L_3$ is isomorphic to
$\bigwedge^3 V$ regardless of $n$, where $V$ is the standard
representation of $\SL_n$.  

In \cite{CF1} we defined a notion of representation stability for each of the sequences of groups $A_n$ given in Table~\ref{table:other}.     We gave some examples, gave some conjectures using this language, and worked out some of the basic theory.    The powerful FI-module point of view was only developed in the special case of $A_n=S_n$.  This is completely missing in general.

\begin{problem}[{\bf FI-theory for other sequences of groups}]
\label{problem:generalfi}
For each of the sequences $A_n=\Sp_{2n}\Z, \GL_n\Z, \SL_n\F_p,\Sp_{2n}\F_p$, work out a theory of 
$\text{FI}_\text{A}$-modules, where:
\begin{enumerate}
\item Finite generation (perhaps with an additional condition) is equivalent to representation stability for $A_n$-represntations, as defined in \cite{CF1}.     
\item The theory gives uniform descriptions (uniform in $n$) of the characters of the examples in Table~\ref{table:other}.  
\item $\text{FI}_\text{A}$-modules satisfy a Noetherian property.
\end{enumerate}
\end{problem}

In \cite{Wi2} J. Wilson extended the theory of FI-modules from the case of $S_n$ to the other two sequences of classical Weyl groups (of type B/C and D); this includes the hyperoctahedral groups.    New phenomena occur here.  For example,  character polynomials must be given with two distinct sets $\{X_i\}, \{Y_i\}$ of variables.   Wilson applies this theory to a number of examples, including the cohomology of the pure string motion groups (see also \cite{Wi1}), the cohomology of various hyperplane arrangements, and 
diagonal co-invariant algebras for Weyl groups.



\begin{thebibliography}{7}

\bibitem{Ar} V.I. Arnol'd, The cohomology ring of the colored
  braid group, {\em Mathematical Notes} 5, no. 2 (1969), 138-140.
 
\bibitem{BM} 
M. Bendersky and J. Miller, Localization and homological stability of configuration spaces,  {\em Quarterly Jour. of Math.}, to appear. arXiv:1212.3596.

\bibitem{BDO}
C. Bowman, M. De Visscher and R. Orellana, The partition algebra and the Kronecker coefficients, preprint, February 2013,  arXiv:1210.5579v6. 

\bibitem{Ca} F. Calegari, The stable homology of congruence subgroups, preprint, November 2013.  arXiv:1311.5190.

\bibitem{Cha}
R. Charney, On the problem of homology stability for
congruence subgroups, {\em Comm. in Alg.}, 12:17-18, 2081-2123 (1984).

\bibitem{Ch}
T. Church, Homological stability for configuration spaces of manifolds, {\em Invent.\  Math.\ } 188 (2012) 2, 465--504.

\bibitem{CE} T. Church and J.\,S. Ellenberg, Homological properties of FI-modules and stability, in preparation.

\bibitem{CEF1} T. Church, J.\,S. Ellenberg and B. Farb, FI-modules: a new approach to stability for $S_n$-representations, preprint, June 2012. arXiv:1204.4533.

\bibitem{CEF2}
T. Church, J.\,S. Ellenberg and B. Farb, Representation stability in cohomology and asymptotics for families of varieties over finite fields, to appear in {\em Algebraic Topology: Applications and New Directions}, AMS Contemp.\ Math.\ Series. 

\bibitem{CEFN}
T. Church, J.\,S. Ellenberg, B. Farb and R. Nagpal, FI-modules over Noetherian rings, {\em Geometry \& Topology}, to appear. 

\bibitem{CF1}
T. Church and B. Farb, Representation theory and homological stability,  {\em Advances in Math.}, 
Vol. 245 (2013), pp. 250--314.

 \bibitem{CF2}
  T. Church and B. Farb, Parameterized Abel--Jacobi maps and abelian cycles in the Torelli group, \emph{Journal of Topology},  5 (2012), no. 1, 15--38. 

\bibitem{CP}
T. Church and A. Putman, Generating the Johnson filtration, preprint, November 2013, arXiv:1311.7150.

\bibitem{CT} F. R. Cohen and L.R. Taylor, Computations of Gel'fand-Fuks cohomology, the cohomology of function spaces, and the cohomology of configuration spaces, in `` Geometric applications of homotopy theory (Proc. Conf., Evanston, Ill., 1977)'', I, pp. 106Ð143, {\em Springer Lect. in Math.}, Vol. 657, 1978.

\bibitem{Co}
R. Cohen, Stability phenomena in the topology of moduli spaces, in ``Surveys in differential geometry,  Vol. XIV, Geometry of Riemann surfaces and their moduli spaces'', pp. 23--56, {\em Surv. Differ. Geom.}, 14, Int. Press, 2009. 

\bibitem{He} D. Hemmer, Stable decompositions for some symmetric group characters arising in braid group cohomology, \emph{J. Combin. Theory Ser. A} 118 (2011), 1136--1139.

\bibitem{JR1} R. Jimenez Rolland, Representation stability for the cohomology of the moduli space ${\cal M}_g^n$, {\em Algebr. Geom. Topol.} 11 (2011), no. 5, 3011--3041. 

\bibitem{JR2} R. Jimenez Rolland, On the cohomology of pure mapping class groups as FI-modules, 
\emph{J. Homotopy Relat. Struct.}, 2013.

\bibitem{Leh} G.I. Lehrer, The $\ell$-adic Cohomology of Hyperplane Complements, \emph{Bull. Lond. Math. Soc.} 24 (1992) 1, 76--82.

\bibitem{LS} G. Lehrer and L. Solomon, On the action of the symmetric group on the cohomology of the complement of its reflecting hyperplanes, \emph{J. Algebra} 104 (1986), no. 2,  410--424.

\bibitem{Li} D.E. Littlewood, Products and plethysms of characters with orthogonal, symplectic and symmetric groups, {\em Canad. J. Math.} 10 1958 17--32. 

\bibitem{Lu} W. L\"{u}ck, \emph{Transformation Groups and Algebraic $K$-Theory}, Lecture Notes in Math., Vol. 1408, Springer-Verlag, Berlin,  1989.

\bibitem{Mac}
  I.G. Macdonald, {\it Symmetric functions and Hall polynomials}, second ed., Oxford Math. Mon., Clarendon Press, Oxford, 1995.
  
\bibitem{Mu} F.D.  Murnaghan, The analysis of the Kronecker product of irreducible representations of the symmetric group, {\em American Jour. of Math.}, Vol. 60, No. 3 (July 1938), pp. 761-784.

\bibitem{Na} R. Nagpal, FI-Modules: Cohomology of modular $S_n$-representations.

\bibitem{Pu}
A. Putman, Stability in the homology of congruence subgroups, arXiv{1201.4876}{4}, revised August 2012.

\bibitem{RW} O. Randal-Williams, Homological stability for unordered configuration spaces, {\em Quarterly Jour. of Math.}, \ Vol. 64, No. 1, pp. 303-326 (2013).
  
\bibitem{SS}
S. Sam and A. Snowden, $\GL$-equivariant modules over polynomial rings in infinitely many variables, arXiv{1206.2233}{2}, revised October 2013.

\bibitem{Sn} A. Snowden, Syzygies of Segre embeddings and $\Delta$-modules, \emph{Duke Math J.} 162 (2013)  2, 225--277. 
  
 \bibitem{Sr}
  B. Srinivasan, \emph{Representations of Finite Chevalley Groups}, Lecture Notes in Math.\ Vol.\ 764, Springer-Verlag, 1979.

\bibitem{To}
B. Totaro, Configuration spaces of algebraic varieties, \emph{Topology} 35 (1996), no. 4, 
1057--1067. 

 \bibitem{VW} R. Vakil and M.M. Wood, Discriminants in the Grothendieck ring, arXiv:1208.3166.

\bibitem{Wi1} J.C.H. Wilson, Representation stability for the cohomology of the pure string motion groups, {\em Alg.\ and Geom.\  Top.} 12 (2012) 909--93.

\bibitem{Wi2} J.C.H. Wilson, $\text{FI}_W$-modules and stability criteria for representations of 
classical Weyl groups, preprint, Sept.\ 2013, arXiv:1309.3817.

\end{thebibliography}
\end{document}